\documentclass[11pt]{article}

\usepackage{latexsym,epsfig,bm}
\usepackage{amssymb}
\usepackage{color}
\usepackage{amsthm}
\usepackage{mathrsfs}

\DeclareSymbolFont{EUR}{U}{eur}{m}{n}
\SetSymbolFont{EUR}{bold}{U}{eur}{b}{n}
\DeclareSymbolFontAlphabet{\eur}{EUR}

\DeclareSymbolFont{EUB}{U}{eur}{b}{n}
\SetSymbolFont{EUB}{bold}{U}{eur}{b}{n}
\DeclareSymbolFontAlphabet{\eub}{EUB}

\DeclareSymbolFont{EUS}{U}{eus}{m}{n}
\SetSymbolFont{EUS}{bold}{U}{eus}{b}{n}
\DeclareSymbolFontAlphabet{\eus}{EUS}

\DeclareSymbolFont{EUF}{U}{euf}{m}{n}
\SetSymbolFont{EUF}{bold}{U}{euf}{b}{n}
\DeclareSymbolFontAlphabet{\euf}{EUF}

\DeclareSymbolFont{AMSb}{U}{msb}{m}{n}
\DeclareSymbolFontAlphabet{\mathbb}{AMSb}

\newcommand{\bfr}{\begin{flushright}}
\newcommand{\efr}{\end{flushright}}
\newcommand{\bfl}{\begin{flushleft}}
\newcommand{\efl}{\end{flushleft}}

\newcommand{\co}{{\rm const}}

\def\Re {{\rm Re\, }}
\def\Im {{\rm Im\,}}

\newcommand{\2}{\hspace{0.2mm}}

\newcommand{\beqn}{\begin{eqnarray}}
\newcommand{\eeqn}{\end{eqnarray}}
\newcommand{\be}{\begin{equation}}
\newcommand{\ee}{\end{equation}}

\newcommand{\E}{{\cal E}}

\newcommand{\bS}{{\bf S}}

\newcommand{\de}{\delta}

\newcommand{\ka}{\kappa}

\newcommand{\om}{\omega}

\newcommand{\toEF}{\stackrel{{\E}\sb F}\lr}
\newcommand{\lr}{\longrightarrow}

\newcommand\p{\partial}
\newcommand{\at}[1]{\vert\sb{\sb{#1}}}

\def\R{{\mathbb R}}

\providecommand\C{{\mathbb C}}
\renewcommand\C{{\mathbb C}}

\newcommand{\abs}[1]{\vert #1 \vert}

\newcommand{\norm}[1]{\Vert #1 \Vert}

\newcommand\sothat{{\rm :}\ }

\newcommand\mod{\mathop{\rm mod}}

\providecommand{\ltor}[1]{
\ifnum #1=1
 {\it i}
\else
\ifnum #1=2
 {\it ii}
\else
\ifnum #1=3
 {\it iii}
\else
\ifnum #1=4
 {\it iv}
\fi
\fi
\fi
\fi
}

\DeclareMathSymbol{\varGamma}{\mathord}{letters}{"00}
\DeclareMathSymbol{\varDelta}{\mathord}{letters}{"01}
\DeclareMathSymbol{\varTheta}{\mathord}{letters}{"02}
\DeclareMathSymbol{\varLambda}{\mathord}{letters}{"03}
\DeclareMathSymbol{\varXi}{\mathord}{letters}{"04}
\DeclareMathSymbol{\varPi}{\mathord}{letters}{"05}
\DeclareMathSymbol{\varSigma}{\mathord}{letters}{"06}
\DeclareMathSymbol{\varUpsilon}{\mathord}{letters}{"07}
\DeclareMathSymbol{\varPhi}{\mathord}{letters}{"08}
\DeclareMathSymbol{\varPsi}{\mathord}{letters}{"09}
\DeclareMathSymbol{\varOmega}{\mathord}{letters}{"0A}

\font\thf cmssdc10 at 11pt
\theoremstyle{plain}
\newtheorem{theorem}{\thf Theorem}[section]

\theoremstyle{definition}
\newtheorem{definition}[theorem]{Definition}

\theoremstyle{remark}
\newtheorem{remark}[theorem]{Remark}

\makeatletter\@addtoreset{equation}{section}
\makeatother

\begin{document}

\title{
On Global Attraction to Solitary Waves for
 the \\Klein-Gordon Equation
Coupled to Nonlinear Oscillator
}

\author{
{\sc Alexander Komech}
\footnote{
On leave from Department of Mechanics and Mathematics,
Moscow State University, Moscow 119899, Russia.
Supported in part
by Max-Planck Institute for Mathematics in the Sciences (Leipzig),
the
Wolfgang Pauli Institute and the Faculty of Mathematics,
Vienna University,
and by DFG Grant (436 RUS 113/615/0-1).}
{\it Faculty of Mathematics, Vienna University,
 Vienna A-1090, Austria}
\medskip\\
{\sc Andrew Komech}
\footnote{
Supported in part
by Max-Planck Institute for Mathematics in the Sciences (Leipzig)
and by the NSF Grant DMS-0434698.
}
{\it
Mathematics Department, Texas A\&M University,
College Station, TX, USA}
}

\date{\today}

\maketitle

\markboth{Alexander Komech, Andrew Komech}
{On Global Attraction to Solitary Waves}
~\\
{\bf Abstract}
The long-time asymptotics is analyzed for all finite energy
solutions to a model $\mathbf{U}(1)$-invariant
nonlinear Klein-Gordon equation in one dimension,
with the nonlinearity concentrated at a point.
Our main result is that
each finite energy solution converges
as $t\to\pm\infty$
to the set of ``nonlinear eigenfunctions''
$\psi(x)e\sp{-i\om t}$.
\medskip\\
{\bf R\'esum\'e.} {\sc Attraction Globale vers des Ondes Solitaires
pour l'\'Equation de Klein-Gordon Coupl\'e \`a
un Oscillateur non Lin\'eaire
}.
On s'int\'eresse aux solutions d'\'energie finie d'une \'equation non
lin\'eaire de Klein-Gordon
$\mathbf{U}(1)$-invariante monodimensionnelle, avec une non lin\'earit\'e
ponctuelle,
et on analyse leur comportement asymptotique aux temps longs.
Le principal r\'esultat que nous avons obtenu est que toute solution
d'\'energie finie
converge pour $t\to\pm\infty$ vers un ensemble de "fonctions propres non
lin\'eaires"
$\psi(x)e\sp{-i\om t}$.
\section{Introduction}
\begin{figure}[htbp]
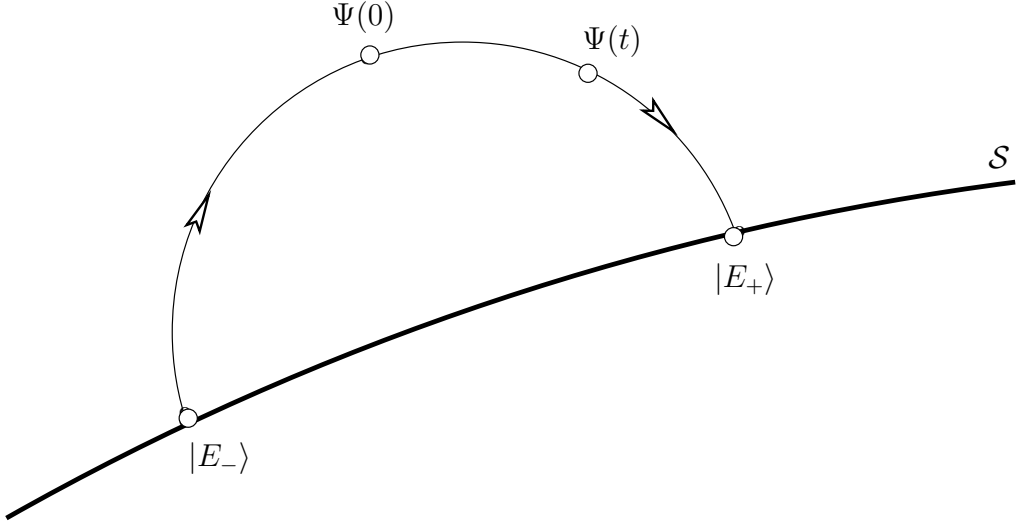

\begin{center}
\input f2.pstex_t
\caption{
Attraction of any trajectory $\Psi(t)$ to the global attractor $\mathcal{S}$.}
\label{fig-gaa}
\end{center}
\end{figure}

We consider the global attractor,
that is, the attracting set for
all finite energy solutions to a model system.
For the first time, we prove that
in a particular $\mathbf{U}(1)$-invariant
dispersive Hamiltonian system
the global attractor is finite-dimensional
and is formed by solitary waves.
The investigation is inspired by Bohr's quantum transitions
(``quantum jumps").
Namely, according to Bohr's postulates,
an unperturbed electron lives forever
in a \emph{quantum stationary state} $\vert E\rangle$
that has a definite value $E$ of the energy.
Under an external perturbation,
the electron
can
jump from one state to another:
$
\vert E\sb{-}\rangle
\longmapsto
\vert E\sb{+}\rangle.
$
The postulate
suggests the dynamical interpretation
of the transitions
as long-time attraction
\begin{equation}\label{ga}
\Psi(t)\lr \vert E\sb\pm\rangle,
\qquad
t\to\pm\infty
\end{equation}
for any trajectory $\Psi(t)$ of the corresponding dynamical system,
where the limiting states
$\vert E\sb\pm\rangle$ generally
depend on the trajectory.
Then the  quantum stationary states
should be viewed as the points
of the \emph{global attractor}
$\mathcal{S}$ which is the set of all
limiting states.
See Figure~\ref{fig-gaa}.
Similar convergence to a global attractor
is well-known for dissipative systems, like Navier-Stokes equations
(see \cite{MR1156492,MR1441312}).
In this context,
the global attractor is formed by the \emph{static stationary states},
and the corresponding asymptotics (\ref{asymptotics})
only holds for $t\to+\infty$ (and with $\om\sb{+}=0$).

Following de Broglie's ideas,
Schr\"odinger identified the stationary states
$\vert E\rangle$
as the solutions
of the
wave equation that have the form
$
\psi(x,t)=\phi\sb\om(x)e\sp{-i\om t}
$, where
$
\om=E/\hbar,
$
and $\hbar$ is Planck's constant.
Then the attraction (\ref{ga})
takes the form of the long-time asymptotics
\begin{equation}\label{asymptotics}
\psi(x,t)
\sim
\psi\sb\pm(x,t)=\phi\sb{\om\sb\pm}(x)e\sp{-i\om\sb{\pm}t},
\qquad
t\to\pm\infty,
\end{equation}
that hold for each finite energy solution.
Our main impetus for considering this problem
was the natural question
whether dispersive Hamiltonian systems
could, in the same spirit,
possess finite dimensional global attractors,
and whether such attractors
are formed by the solitary waves.
We prove such a global attraction
for a model nonlinear Klein-Gordon equation
\begin{equation}\label{KG-0}
\ddot\psi(x,t)
=\psi''(x,t)-m^2\psi(x,t)+\de(x)F(\psi(0,t)),
\qquad
\quad x\in\R.
\end{equation}
Here $m>0$, $\psi(x,t)$ is a continuous complex-valued wave function,
and $F$ is a nonlinearity.
The dots stand for the derivatives in $t$,
and the primes for the derivatives in $x$.
We assume that equation (\ref{KG-0}) is $\mathbf{U}(1)$-invariant;
that is,
$
F(e\sp{i\theta}\psi)=e\sp{i\theta}F(\psi),
$
$
\theta\in\R.
$
\medskip\\
Let $\mathcal{S}$ be the set of all functions
$\phi\sb\om(x)\in H\sp 1(\R)$
with $\omega\in\C$,
so that $\phi\sb\om(x)e\sp{-i\omega t}$
is a solution to (\ref{KG-0}).
Our main result is the following long-time asymptotics
(cf. (\ref{asymptotics})) for {\it nonlinear polynomial functions} $F(\psi)$:
\begin{equation}\label{attraction}
\psi(\cdot,t)\longrightarrow\mathcal{S},
\qquad t\to\pm\infty,
\end{equation}
where the convergence holds in local energy seminorms.
In the linear case, when $F(\psi)=a\psi$ with $a\in\R$,
there is no global attraction
to $\mathcal{S}$ if $a>0$, although
the attraction holds if $a\le 0$
(see Remark~\ref{remark-linear}).
Although we proved the attraction (\ref{attraction}) to $\mathcal{S}$,
we have not proved the attraction to a particular point of $\mathcal{S}$,
falling short of proving (\ref{asymptotics}).
Hypothetically, a solution can be drifting along $\mathcal{S}$,
keeping asymptotically close to it,
but never stopping at a single point of $\mathcal{S}$.
Let us comment on related earlier results:

\noindent {\it i})
The asymptotics of type
(\ref{ga})
with $\vert E\sb{\pm}\rangle=0$
were discovered in the scattering theory
\cite{MR0153967, MR0233062,MR0303097,MR654553,MR824083}.
In this case, the attractor $\mathcal{S}$ consists of the zero solution only,
and the asymptotics mean well-known \emph{local energy decay}.

\noindent {\it ii})
The \emph{global attraction}
(\ref{ga})
with
$\vert E\sb\pm\rangle\ne 0$
was established first in
\cite{MR1203302e,MR1434147,MR1726676,MR1748357}
for a number of nonlinear
wave problems.
There the attractor $\mathcal{S}$
is the set of all \emph{static} stationary states.
Let us mention that this set could be infinite
and contain continuous components.

\noindent {\it iii})
First results on the asymptotics of type (\ref{asymptotics}),
with $\omega\sb\pm\ne 0$
were obtained for nonlinear $\mathbf{U}(1)$-invariant
Schr\"odinger equations
in the context of asymptotic stability.
This establishes asymptotics of type (\ref{asymptotics})
but only for solutions close to the solitary waves,
proving the existence of a \emph{local attractor}.
This was first done in
\cite{MR1071238},
and then developed in
\cite{MR1199635,MR1681113,MR1835384} and others.

\medskip

Let us mention that the global
attraction (\ref{asymptotics})
for equation (\ref{KG-0}) with $m=0$
follows from
\cite{MR1203302e}; In that case, $\omega\sb\pm=0$.
Our proofs for the case $m>0$
are quite different from the approach used in \cite{MR1203302e},
and are based
on a nonlinear spectral analysis
of \emph{omega-limit trajectories} at $t\to\pm\infty$.


\section{Main Results}
\setcounter{equation}{0}
We consider the Cauchy problem for the
equation (\ref{KG-0}).
We define
$\Psi(t)=
\left[\!\scriptsize{\begin{array}{c}
\psi(x,t)\\\pi(x,t)\end{array}}\!\right]$
and write the Cauchy problem
in the vector form:
\begin{equation}\label{KG-cp}
\dot\Psi(t)
=
\left[\begin{array}{cc}0&1\\\p\sb x\sp 2-m^2&0\end{array}\right]
\Psi(t)
+
\de(x)\left[\begin{array}{c}0\\F(\psi)\end{array}\right],
\qquad
\Psi\at{t=0}
=\Psi\sb 0
\equiv\left[\begin{array}{c}\psi\sb 0\\\pi\sb 0\end{array}\right].
\end{equation}

\begin{definition}
\begin{enumerate}
\item
${\E}$ is the Hilbert space of the states
$\Psi
=(\psi(x),\pi(x))$,
with the norm
\[
\norm{\Psi}\sb{\E}^2
:=
\norm{\psi'}\sb{L\sp 2}^2+\norm{\psi}\sb{L\sp 2}^2+\norm{ \pi}\sb{L\sp 2}^2,
\qquad
{\rm where}
\ \ L\sp 2=L\sp 2(\R).
\]
\item
${\E}\sb F$ is the space ${\E}$ endowed
with the Fr\'echet topology defined by the seminorms
\[
\norm{\Psi }\sb{\E,R}^2
:=
\norm{\psi'}\sb{L\sp 2\sb R}^2+\norm{\psi}\sb{L\sp 2\sb R}^2+
\norm{ \pi}\sb{L\sp 2\sb R}^2,
\qquad
{\rm where}
\ \ L\sp 2\sb R=L\sp 2(-R,R),\ \ R>0.
\]
\end{enumerate}
\end{definition}

We assume that
the oscillator force $F$ admits a real-valued
potential:
$
F(\psi)=-\nabla U(\psi),\quad\psi\in\C
$,
where
$
U\in C\sp 2(\C)
$, and
the gradient is taken with respect to $\Re\psi$ and $\Im\psi$.
Then equation (\ref{KG-cp})
formally can be written as a Hamiltonian system.
We assume that the potential $U(\psi)$ is $\mathbf{U}(1)$-invariant,
where $\mathbf{U}(1)$ stands for the unitary group
$e\sp{i\theta}$, $\theta\in\R\mod 2\pi$:
Namely, we assume that
there exists $u\in C\sp 2(\R)$ such that
$U(\psi)=u(\abs{\psi}\sp 2)$,
$\psi\in\C$.

\begin{theorem}\label{theorem-well-posedness}
Assume that
$U(\psi)\ge A-B\abs{\psi}^2$,
where $A,\,B\in\R$ and $B<m$.
Then
for every $\Psi\sb 0\in {\E}$  the Cauchy problem
(\ref{KG-cp}) has a unique solution
$\Psi(t)=(\psi(x,t),\pi(x,t))\in C(\R,{\E})$.
The energy is conserved:
$$
\frac 1 2
\int\sb\R
\left(
|\pi(x,t)|\sp 2+|\psi'(x,t)|^ 2+m^2 |\psi(x,t)|\sp 2
\right)
\,dx
+U(\psi(0,t))
=\co,~~~~~~t\in\R,
$$
and
a priori bound
$\norm{\Psi(t)}\sb{\mathcal{E}}\le C(\norm{\Psi\sb 0}\sb{\mathcal{E}})$
holds for
$t\in\R$.

\end{theorem}

\begin{definition}\label{def-solitary-waves}
\begin{enumerate}
\item
The solitary waves of equation (\ref{KG-cp})
are solutions of the form
\begin{equation}\label{solitary-waves}
\Psi(t)=\Phi\sb\omega e\sp{-i\omega t},
\qquad
{\rm where }\quad
\omega\in\C,
\quad
\Phi\sb\omega
=
\left[\begin{array}{r}\phi\sb\om(x)
\\
-i\omega\phi\sb\om(x)\end{array}
\right],
\quad
\phi\sb\omega\in H\sp 1(\R).
\end{equation}
\item
The solitary manifold
is the set
$
\bS=
\left\{
\Phi\sb\omega\sothat\omega\in\C
\right\}
$
of all amplitudes $\Phi\sb\omega$.
\end{enumerate}
\end{definition}

The profiles of the solitary waves have the form
$\phi\sb\omega(x)=C e^{-\ka\abs{x}}$,
where
$C\in\C$, $\kappa\ge 0$, and $\omega\in\C$
are related by
the \emph{linear dispersion relation} $\omega^2=m^2-\kappa^2$
and the \emph{coupling identity}
$2\kappa C=F(C)$.
Thus, $\bS$
is generically a two-dimensional real submanifold of $\E$
that can be parametrized by the corresponding
complex amplitudes $C$.

\begin{theorem}
\label{main-theorem}
Let the nonlinearity $F(\psi)$ satisfy
$F(\psi)=-\nabla U(\psi)$,
where
\begin{equation}\label{a-u}
U(\psi)
=\sum\limits\sb{n=0}\sp{N}u\sb n\abs{\psi}\sp{2n},
\quad
N\ge 2;
\quad
u\sb n\in\R,
\quad
u\sb N > 0.
\end{equation}
Then for any $\Psi\sb 0\in \E$
the solution $\Psi(t)\in C(\R,{\E})$
to the Cauchy problem {\rm (\ref{KG-cp})}
with $\Psi(0)=\Psi\sb 0$
converges to the set $\bS$ in the space $\E\sb{F}$:
\begin{equation}\label{cal-A}
\Psi(t)\toEF \bS,
\quad t\to \pm\infty.
\end{equation}
\end{theorem}

The assumption (\ref{a-u})
that the nonlinearity is polynomial
is crucial in our argument:
It will allow to apply the Titchmarsh convolution theorem.
Under this assumption,
nonzero solitary waves (\ref{solitary-waves})
correspond only to real values of $\om\in(-m,m)$.
\begin{remark}\label{remark-linear}
In the linear case, when $F(\psi)=a\psi$ and $a>0$,
the equation admits two linearly independent
solutions
$\psi_\pm(x,t)=e^{-a\abs{x}/2}
e^{-i\2 \om_\pm t}$
with $\om_\pm=\pm\sqrt{m^2-a^2/4}$ if $m\ne a/2$,
and $e^{-m\abs{x}}$, $te^{-m\abs{x}}$ if $m= a/2$.
Hence the global attraction (\ref{cal-A}) fails
because of the superposition principle.
For $a\le 0$ we have $\bS=\{0\}$,
and the attraction (\ref{cal-A}) holds.

\end{remark}

\noindent
{\bf Strategy of the proof of Theorem \ref{main-theorem}}
For the Klein-Gordon equation with $m>0$,
the dispersive relation $\omega\sp 2=k\sp 2+m^2$
results in the group velocities $v=\omega'(k)=k/\sqrt{k^2+m^2}$,
so every velocity $0\le\abs{v}<1$ is possible.
This complicates considerably the investigation of the energy propagation,
so
the approach of \cite{MR1203302e}
built on the fact that the group velocity was $\abs{v}=1$
no longer works.

We prove the absolute continuity
of the spectrum of the solution  for $\abs{\omega}>m$.
This observation is similar to
the well-known Kato Theorem.
The proof is not obvious
and relies on the complex Fourier-Laplace transform and
the Wiener-Paley arguments.

We then split the solution into two components: Dispersive and bound,
with the frequencies  $\abs{\omega}>m$ and $\omega\in [-m,m]$,
respectively.
The dispersive component is an oscillatory integral of plane waves,
while the bound component
is a superposition of exponentially decaying functions.
The stationary phase argument leads to a local decay of
the dispersive component,
due to the absolute continuity of its spectrum.
This reduces the long-time behavior of the solution to the
behavior of the bound component.

Next, we establish the spectral representation
for the bound component. For this, we need to know
an optimal regularity of the corresponding spectral measure;
We have found out
that the spectral measure belongs to the space of
\emph{quasimeasures}
 which are Fourier transforms of bounded continuous functions.
The spectral representation
implies compactness in the space of quasimeasures,
which in turn leads to the existence
of \emph{omega-limit trajectories}  at $t\to\pm\infty$.

Further,  we prove that
an omega-limit trajectory itself
satisfies the nonlinear equation (\ref{KG-0}),
and this implies
the crucial spectral inclusion:
The spectrum of the nonlinear term
is included in the spectrum of the omega-limit trajectory.
We then reduce the spectrum of this limiting trajectory
to a single harmonic $\omega\sb{\pm}\in[-m,m]$
using the Titchmarsh convolution theorem
\cite{titchmarsh} (see also \cite[Theorem 4.3.3]{MR1065136}).
In turn, this means that any omega-limit trajectory
lies in the manifold $\bS$ of the solitary waves,
proving that $\bS$ is the global attractor.

\bibliographystyle{ubk}
\bibliography{ubk-mathsci,ubk-local}

\end{document}